\documentclass[12pt]{article}
\usepackage{amsmath,amssymb,amsthm,graphicx,float}
\usepackage[]{algorithm2e}

\theoremstyle{plain}
\newtheorem{theorem}{Theorem}

\newtheorem{lemma}[theorem]{Lemma}

\newtheorem{definition}[theorem]{Definition}

\newcommand\numberthis{\addtocounter{equation}{1}\tag{\theequation}}

\newcommand{\zp}{\mathbb{F}_p}
\DeclareMathOperator{\Ch}{Ch}

\newcommand{\oo}{\mathcal{O}}

\title{Arithmetic Progressions of Length Three in Multiplicative Subgroups of $\zp$}
\author{Jeremy F.~Alm}
\date{January 2018}

\begin{document}

\maketitle

\begin{abstract}

In this paper, we give an algorithm for detecting non-trivial 3-APs in multiplicative subgroups of $\mathbb{F}_p^\times$ that is substantially more efficient than the naive approach.  It follows that certain Var der Waerden-like numbers  can be computed in polynomial time.  
\end{abstract}

\section{Introduction}

Additive structures inside multiplicative subgroups of $\zp^\times$ have recently received attention.  Alon and Bourgain \cite{MR3213827} study solutions to $x+y=z$ in $H<\zp^\times$, and Chang \cite{MR3722263} studies arithmetic progressions  in $H<\zp^\times$. In this paper, we define a Van der Waerden-like number for $H<\zp^\times$ of index $n$, and give a polynomial-time algorithm for determining such numbers.

\begin{definition}
Let $VW_3^\times(n)$ denote the least prime $q\equiv 1 \pmod{n}$ such that for all primes $p\equiv 1 \pmod{n}$ with $p\geq q$, the multiplicative subgroup of $\zp^\times$ of index $n$ contains a mod-$p$ arithmetic progression of length three.
\end{definition}

Our main results are the following two theorems:
\begin{theorem}\label{th1}
$VW_3^\times(n) \leq (1+\varepsilon )n^4$ for all sufficiently large $n$ (depending on $\varepsilon$). In particular, $VW_3^\times(n) \leq 1.001n^4$ for all $n\geq 45$. 
\end{theorem}

\begin{theorem}\label{thm:alg}
$VW_3^\times(n)$ can be determined by an algorithm that runs in $\oo(\frac{n^8}{\log n})$ time.
\end{theorem}

Chang \cite{MR3722263} proves that if $H<\zp^\times$ and $|H|> cp^{3/4}$, then $H$ contains non-trivial 3-progressions. This implies our Theorem \ref{th1} with $(1+\varepsilon)n^4$ replaced by $cn^4$.  We prove our Theorem \ref{th1} because we need to make the constant explicit.

\section{Proof of Theorem \ref{th1}}
\begin{proof}

We use one of the basic ideas of the proof of Roth's Theorem on 3-progressions \cite{Roth}. Let $A\subseteq \zp$ with $|A|=\delta p$. Note that a 3-progression   is a solution inside $A$ to the equation $x+y=2z$.  Let $\mathcal{N}$ be the number of (possibly trivial) solutions to $x+y=2z$ inside $A$.  We have that 
 \begin{equation} \label{eq1}
  \frac{1}{p}\sum^{p-1}_{k=0} e^{\frac{-2\pi ik}{p}x}= 
     \begin{cases}
     1, & \text{ if } x\equiv 0 \pmod p ;\\
     0, &\text{ if } x\not\equiv 0 \pmod p.
 \end{cases}
 \end{equation}

 Because of \eqref{eq1}, we have
 
 \begin{equation}\label{eq2}
 \mathcal{N}=\sum_{x\in A}\sum_{y\in A}\sum_{z\in A}\frac{1}{p}\sum_{k=0}^{p-1} e^{\frac{-2\pi ik}{p}(x+y-2z)}
 \end{equation}
 
 Rearranging \eqref{eq2}, we get

 \begin{align*}
 &\phantom{=} \ \ \frac{1}{p}\sum_{k=0}^{p-1}\sum_{x\in A}\sum_{y\in A}\sum_{z\in A} e^{\frac{-2\pi ik}{p}x}\cdot e^{\frac{-2\pi ik}{p}y}\cdot e^{\frac{2\pi ik}{p}z}\\
 &= \frac{1}{p}\sum_{k=0}^{p-1}\left[  \sum_{x\in A}e^{\frac{-2\pi ik}{p}x} \ \cdot \ \sum_{y\in A}e^{\frac{-2\pi ik}{p}y} \ \cdot \  \sum_{z\in A}e^{\frac{2\pi ik}{p}2z}\right]\\
 &= \frac{1}{p}\sum_{k=0}^{p-1} \left[\sum_{x\in\zp} \Ch_A(x) e^{\frac{-2\pi ik}{p}x} \ \cdot \  \sum_{y\in\zp} \Ch_A(y) e^{\frac{-2\pi ik}{p}y} \ \cdot \ \sum_{z\in\zp} \Ch_A(-2z) e^{\frac{2\pi ik}{p}z}\right]\\
  &= \frac{1}{p}\sum_{k=0}^{p-1} \hat{\Ch}_A(k)^2\cdot\hat{\Ch}_A(-2k), \numberthis\label{eq3}\\
 \end{align*}
 where $\Ch_A$ denotes the characteristic function of $A$, and $\hat{f}$ denotes the Fourier transform of $f$, 
 \[
   \hat{f}(x) = \sum^{p-1}_{k=0}f(k)  e^{\frac{-2\pi ik}{p}x} .
 \]
 Now we can pull out the $k=0$ term from \eqref{eq3}:
 
 \begin{align*}
 \eqref{eq3} &= \frac{1}{p}\hat{\Ch}(0)^3+\frac{1}{p}\sum_{k=1}^{p-1}\hat{\Ch}_A(k)^2\cdot\hat{\Ch}_A(-2k) \\
 &= \frac{|A|^3}{p} + \frac{1}{p}\sum_{k=1}^{p-1}\hat{\Ch}_A(k)^2\cdot\hat{\Ch}_A(-2k)\\
 &=\delta^3p^2 + \frac{1}{p}\sum_{k=1}^{p-1}\hat{\Ch}_A(k)^2\cdot\hat{\Ch}_A(-2k).
 \end{align*}
 
 Let's call $\delta^3p^2$ the \emph{main term}, and $\frac{1}{p}\sum_{k=1}^{p-1}\hat{\Ch}_A(k)^2\cdot\hat{\Ch}_A(-k)$ the \emph{error term}.  We now bound this error term.
 
 Suppose $0<\alpha <1$ and $|\hat{\Ch}_A(k)|\leq\alpha p$ for all $0 \neq k\in\zp$. In this case, we say that $A$ is $\alpha$-\emph{uniform}. Then
 \begin{align*}
 \left|\frac{1}{p}\sum_{k=1}^{p-1}\hat{\Ch}_A(k)^2\cdot\hat{\Ch}_A(-2k)\right| &\leq \frac{1}{p}\max |\hat{\Ch}_A(k)|\cdot  \left|\sum_{k=1}^{p-1}\hat{\Ch}_A(k)^2\right|\\
 &\leq \alpha  \left|\sum_{k=1}^{p-1}\hat{\Ch}_A(k)^2\right|\\
 &\leq \alpha p\left|\sum_{k=1}^{p-1}\Ch_A(k)^2\right|\\
 &\leq \alpha\delta p^2.
 \end{align*}
 
 Therefore $\mathcal{N}\geq\delta^3 p^2-\alpha\delta p^2$. Subtracting off the trivial solutions gives $\mathcal{N}-\delta p\geq\delta^3 p^2-\delta p - \alpha\delta p^2$. Hence there is at least one non-trivial solution if 
 
 \[
 \delta^3 p^2 > \delta p + \alpha\delta p^2.
 \]

 Let $A = H$ be a multiplicative subgroup of $\zp$ of index $n$. As is well-known (see for example \cite[Corollary 2.5]{Schoen}), if $H$ is a multiplicative subgroup of $\zp^\times$, then $H$ is $\alpha$-uniform for $\alpha \leq p^{-1/2}$. Thus it suffices to have
 \begin{align}
     \delta^3 p^2\geq \delta p + p^{-1/2}\delta p^2 &\Longleftrightarrow \delta^3 p^2\geq \delta p + \delta p^{3/2}\\
     &\Longleftrightarrow \delta^2 p\geq 1 +  p^{1/2}\\
     &\Longleftrightarrow (p-1)^2\geq n^2 p (1 +  p^{1/2}) \label{star}
 \end{align}
 where the last line follows from $\delta = (p-1)/(np)$. It is straightforward to check that \eqref{star} is satisfied by $p=(1+\varepsilon)n^4$ for sufficiently large $n$.

\end{proof}

\begin{figure}[H]
    \centering
    \includegraphics{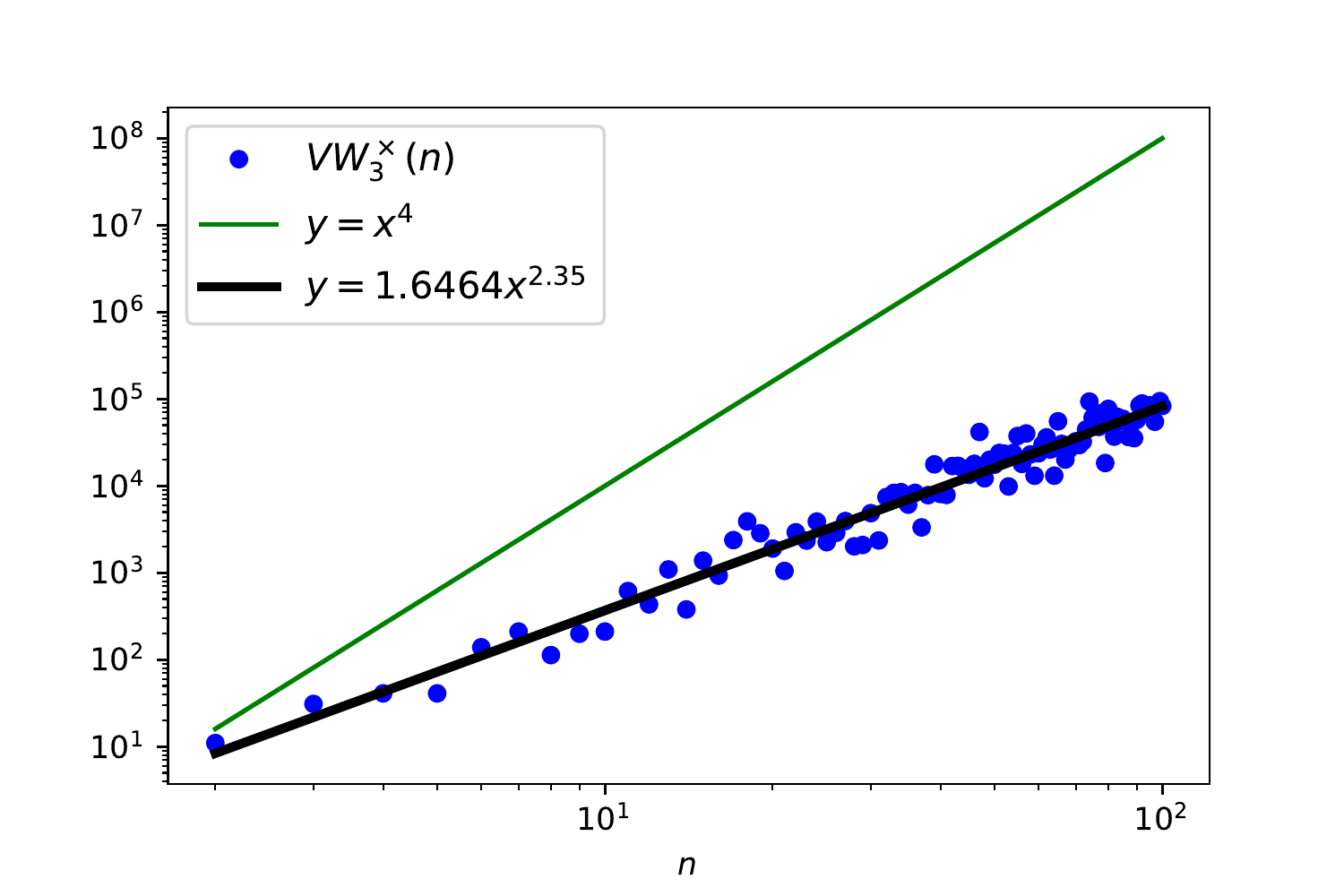}
    \caption{$VW^\times_3(n)$ for $n\leq 100$}
    \label{fig:data}
\end{figure}

The data gathered for $VW^\times_3(n)$, $n\leq 100$, suggest that the exponent of 4 on $n$ is too large; see Figure \ref{fig:data}. These data are available at \texttt{www.oeis.org}, sequence number A298566.

\section{A More General Framework}
Before we establish our algorithm, it will helpful to generalize to arbitrary linear equations in three variables over $\zp$.  Suppose we're looking for solutions to $ax + by = cz$
    in $H<\zp^\times$, for fixed $a,b,c\in \zp^\times$. There is a solution just in case $(aH + bH) \cap cH $ is nonempty.  
    
    The following result affords an algorithmic speedup in counting solutions to $ax + by = cz$ inside $H$:
    
    \begin{lemma}\label{lem}
    For $a$, $b$, $c\in\zp^\times$ and  $H < \zp^\times$,
    \[
        (aH+bH)\cap cH \neq \varnothing  \text{ if and only if } (c - aH)\cap bH \neq \varnothing.
    \]
    \end{lemma}
    
    Notice that while the implied computation on the left side of the biconditional is $\oo(p^2)$, the one on the right is $\oo(p)$, since we compute $|H|$ subtractions and $|H|$ comparisons. (We consider the index $n$ fixed.) 
    
    \begin{proof} 
    Let $H=\{ g^{kn} : 0 \leq k < (p-1)/n \}$, where $n$ is the index of $H$ and $g$ is a primitive root modulo $p$. Fix $a,b,c \in\zp$. 
    
    For the forward direction, suppose $ (aH+bH)\cap cH\neq\varnothing$, so there are $x,y,z \in H$ such that $ax + by = cz$. Then $by = cz-ax$. Multiplying by $z^{-1}\in H$ yields $b(yz^{-1}) = c - a(xz^{-1})$. Therefore $(c - aH)\cap bH \neq \varnothing$.  The other direction is similar.
    \end{proof}
    
    Lemma \ref{lem} allows us to detect solutions to linear equations in linear time. The caveat for the case $a=b=1$, $c=2$ is that $H+H$ \emph{always} contains $2H$, since $h+h=2h$ for all $h\in H$; these solutions correspond to the trivial 3-APs $h,h,h$. (Similarly, $(2 - H)\cap H$ is always nonempty, since $1\in H$ and $2-1=1$.)  To account for this, we simply consider $H'$ = $H\setminus \{1\}$, and calculate $(2 - H')\cap H'$ instead. 
    
\section{Proof of Theorem \ref{thm:alg}}
  
  \begin{proof}
  Here is the algorithm.
  
\begin{algorithm}[H]\label{alg}
 \KwData{An integer $n>1$}
 \KwResult{ The value of $VW_3^\times(n)$}
 
 Let $\mathcal{P} = \{ p \text{ prime} : p \leq (1+\varepsilon)n^4,\ p \equiv 1 \pmod{n} \}$.
 
 Set $p_0 = 1$.
 
 Set Prev\_boolean = False and Current\_boolean = True. 
 
 \For{$p \in \mathcal{P}$}{
    Let $H$ be the subgroup of $\zp^\times$ of index $n$.
    
    Set Current\_boolean to True if $(2 - H')\cap H'$ is non-empty, and False otherwise.

        \If{Current\_boolean is True and Prev\_boolean is False}{
            set $p_0 = p$. 
            }
        Set Prev\_boolean to the value of Current\_boolean.

 }
 
 Return $p_0$\\
 
\caption{Algorithm for determining $VW^\times_3(n)$ }
\end{algorithm}

We now argue that Algorithm \ref{alg} runs in $\oo\left(\frac{n^8}{\log n}\right)$ time. Since calculating $(2 - H')\cap H'$ is $\oo(p)$ for each prime $p$, our runtime is bounded by 

\[
    \sum_{\substack{p \leq (1+\varepsilon)n^4 \\ p \equiv 1 \pmod{n}}} \oo(p) = \oo \left( \sum_{\substack{p \leq (1+\varepsilon)n^4 \\ p \equiv 1 \pmod{n}}} p  \right).
\]

A standard estimate on the prime sum 
\[
        \sum_{\substack{p \leq x \\ p \equiv 1 \pmod{n}}} p
\]
 is asymptotically $\frac{x^2}{\varphi(n)\log x}$, giving 
 
 \begin{align*}
      \oo \left( \sum_{\substack{p \leq (1+\varepsilon)n^4 \\ p \equiv 1 \pmod{n}}} p  \right) 
                                            &= \oo\left(\frac{n^8}{\varphi(n)\log(n^4)}\right) \\
                                            &= \oo\left(\frac{n^8}{\log(n)}\right) \\
 \end{align*}
 as desired.
  \end{proof}  
  
\begin{figure}[H]
    \centering
    \includegraphics{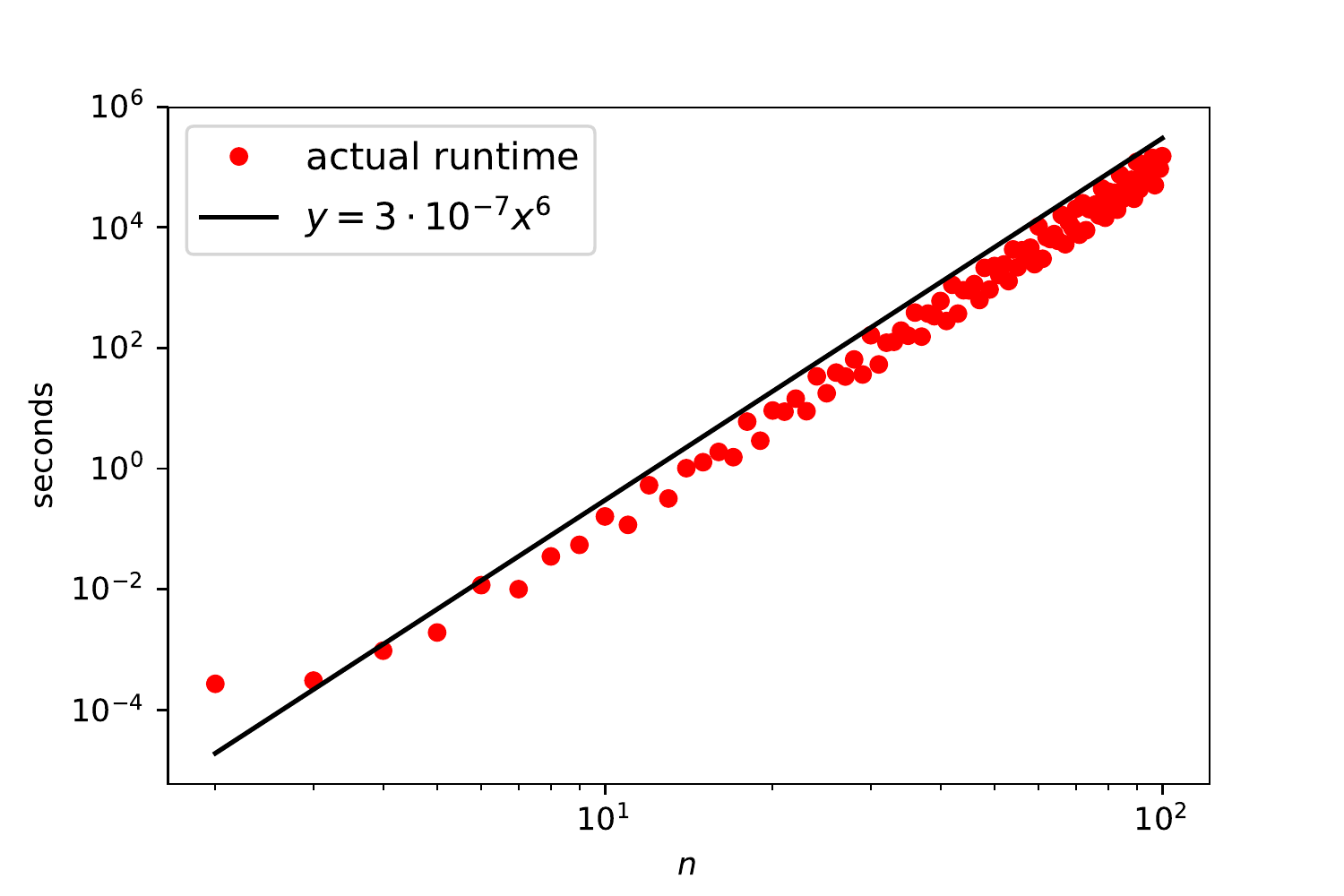}
    \caption{Runtime in seconds to determine $VW_3^\times(n)$}
    \label{fig:timing}
\end{figure}

  Our timing data suggest that the correct runtime might be more like $\oo(n^6)$; see Figure \ref{fig:timing}.

\section{Further Directions}

For any $a,b,c\in \mathbb{Z}^+$, we can define an analog to $VW^\times_3(n)$ by considering the equation $ax+by=cz$ instead of $x+y=2z$.  (Assume $p$ is greater than $a$, $b$, and $c$.) The bound from Theorem \ref{th1} stays the same if $a+b=c$ and goes down to $n^4+5$ otherwise.  But as suggested by the data in Figure \ref{fig:data}, these bounds are not tight.  How does the choice of $a$, $b$, and $c$ affect the growth rate of the corresponding Van der Waerden-like number? Clearly $VW^\times_3(n)$ is not monotonic, but it appears to bounce above and below some ``average'' polynomial growth rate. Will that growth rate vary with the choice of  $a$, $b$, and $c$? Does it depend on whether $a+b=c$ only?

\section{Acknowledgements}

The author wishes to thank Andrew Shallue and Valentin Andreev for many productive conversations.

% \bibliographystyle{plain}
% \bibliography{refs.bib}

\end{document}